\documentclass[a4paper,12pt]{article}
\usepackage{amssymb}
\usepackage{amsmath}

\title{Resolution of Veronese Embedding of plane curves}
\author{Aaloka Kanhere \\ {\small aalokakanhere@gmail.com}\\ Homi Bhabha Center for Science Education}
\begin{document}
\maketitle

{\footnotesize \textbf{Abstract:}
Let $C$ be a smooth (irreducible) curve of degree $d$ in $\mathbb{P}^{2}$. Let $\mathbb{P}^{2} \hookrightarrow \mathbb{P}^{5}$ be the Veronese embedding and let $\mathcal{I}_{C}$ denote the homogeneous ideal of $C$ on $\mathbb{P}^{5}$. In this note we explicitly write down the minimal free resolution of $\mathcal{I}_{C}$ for $d\geq 2$.
}\\
\begin{center}
\textbf{1. Introduction}
\end{center}

In \cite{L}, the author has remarked, " It is very exceptional to be able to construct the whole resolution explicitly, let alone to be able to do so by hand!." This remark of Lazarfeld motivated us to try to explicitly calculate whole resolutions of projective varieties.\\

In this paper I have explicitly calculated the whole resolutions of the Veronese embedding of plane curves. I look at the even and odd degree curves separately and get the explicit resolution for both.\\

Let $C$ be a smooth and irreducible projective curve and $L$ be an ample line
bundle on $C$, generated by its global sections. Then $L$ determines a morphism
\[ \Phi_{L}  :  C \longrightarrow \mathbb{P}(H^{0}(L))=\mathbb{P}^{r}\]
where $r = h^{0}(L)-1$. Also we have that if $L$ is very ample, then $\Phi_{L}$
is an embedding. \\

Let $\mathcal{I}_{C}$ be the homogeneous ideal of $C$ in $\mathbb{P}^{r}$
and $S$, homogeneous coordinate ring of the projective space, $\mathbb{P}^{r}$\\
Let $R= S/\mathcal{I}_{C}$, then the minimal graded free resolution of $R$
is the following exact sequence of free modules:
\[0 \rightarrow E_{n} \stackrel{\alpha_{n}}{\rightarrow} \ldots \stackrel
{\alpha_{3}} {\rightarrow} E_{2} \stackrel{\alpha_{2}}{\rightarrow} E_{1}
\stackrel{\alpha_{1}}{\rightarrow} E_{0} \rightarrow R \rightarrow 0  \ldots
(A)\]

where each $E_{i}$ is a direct sum of twists of $S$, i.e.
\[ E_{i} = \oplus_{j} S(-a_{ij}) , \]

And the maps, $\alpha_{i}$'s in the above exact sequence are given by matrices of homogeneous forms and none of the entries in the above matrices are non-zero constants. Note that $E_{0} = S$ and the image of $\alpha_{1}$ is the ideal of $S$, $\mathcal{I}_{C}$. \\

In this note we look at $C$, a smooth(irreducible) curve of degree $d$ such
that $C \hookrightarrow \mathbb{P}^{2}$(here $L$ is $\mathcal{O}_{C}(1)$). Now because of the Veronese embedding,
we get an embedding of $C$ in $\mathbb{P}^{5}$ which is nothing but the
embedding of $C$ in $\mathbb{P}^{5}$ due to the very ample line bundle,
$\mathcal{O}_{C}(2)$. We explicitly calculate minimal free resolution of $\mathcal{I}_{C}$ and in particular get the equations defining $C$ in $\mathbb{P}^{5}$.\\
Most of the definitions in this note are from \cite{A} and \cite{H}.\\
\\
\\
\begin{center}
\textbf{Notations}
\end{center}

The \emph{first syzygy module} is defined as the image of $\alpha_{2}$ in $E_{1}$ in the exact sequence (A) and is denoted by ${\textrm{Syz}}^{1} (\mathcal{I}_{C})$. \\

The $\mathit {k^{th}} $ \emph{syzygy module} is defined inductively to be the module of syzygies of the $ \mathit{ (k-1)^{st} }$ \emph{syzygy module}. Hence we have the following inductive relation:
\[\mathit{ {\textrm{Syz}}^{k}(\mathcal{I}_{C}) = {\textrm{Syz}}^{1}({\textrm{Syz}}^{k-1}(\mathcal{I}_{C})) } \]
\\

\begin{center}
\textbf{2. Resolutions of Veronese Embedding}\\
\end{center}

Consider $\sigma: \mathbb{P}^{2} \rightarrow \mathbb{P}^{5}$ such that for $p=(a_{0},a_{1},a_{2}) \in \mathbb{P}^{2}$,
\[\sigma(p)= (a_{0}^{2},a_{0}a_{1},a_{0}a_{2},a_{1}^{2},a_{1}a_{2},a_{2}^{2} )\]

This is called the \emph{Veronese embedding} of $\mathbb{P}^{2}$ in $\mathbb{P}^{5}$ \cite{H}.\\

Now if $x_{00},x_{01},x_{02},x_{11},x_{12},x_{22}$ denote homogeneous coordinates on
$\mathbb{P}^{5}$, then one has a description of $\sigma(\mathbb{P}^{2})$ as the
zeros of the six minors of the following  $3 \times 3$ symmetric matrix.\\
\[ \left( \begin{array}{cccc}
           x_{00} &x_{01} &x_{02}\\
            x_{01} &x_{11} &x_{12}\\
              x_{02} &x_{12} &x_{22}\\
              \end{array} \right) \]
\\
Moreover we also get a map,
\[\theta :  k[x_{00},x_{01},x_{02},x_{11},x_{12},x_{22}] \rightarrow
 k[x_{0},x_{1},x_{2}] \] such that,
$\theta(x_{ij}) = x_{i}x_{j}$  $\forall 0 \leq i \leq j \leq 2 $.
\\
Also the defining equations of this embedding are:
\[ \begin{array}{l}
    \Delta_{00} = x_{11}x_{22}-x_{12}^{2} \\
    \Delta_{01} = x_{01}x_{22}-x_{12}x_{02}\\
    \Delta_{02} = x_{01}x_{12}-x_{02}x_{11}\\
    \Delta_{11} = x_{00}x_{22}-x_{02}^{2} \\
    \Delta_{12} = x_{00}x_{12}-x_{02}x_{01}\\
    \Delta_{22} = x_{00}x_{11}-x_{01}^{2}\\
                      \end{array}\]
 Notice that,\\
ker($\theta$) = $<  \hspace{0.1in} \Delta_{i,j}$ , $\forall 0 \leq i \leq j \leq 2 \hspace {0.1in}>$
\\
\\

From now we will denote $k[x_{00},x_{01},x_{02},x_{11},x_{12},x_{22}] $ as $S$. And for $d \in \mathbb{Z}$, $S(d)$ is the graded $S$ module such that $S(d)_{n} = S_{d+n}$
\\
\textbf{Theorem 1} : \cite{O} The ideal $\mathcal{I}_{\mathbb{P}^{2}}$ of
$\sigma(\mathbb{P}^{2})$ in $\mathbb{P}^{5}$ has the following resolution.
\begin{equation}
0 \rightarrow S(-4)^{\oplus 3}
\stackrel{M_{3}}{\rightarrow}
S(-3)^{\oplus 8} \stackrel{M_{2}} {\rightarrow}
S(-2)^{\oplus 6} \stackrel{M_{1}}{\rightarrow}
\mathcal{I}_{\mathbb{P}^{2}} \rightarrow 0
\end{equation}
where,\\
\[
M_{1} = \left[ \begin{array}{rrrrrr} \Delta_{00} ,& \Delta_{01},& \Delta_{02},
,& \Delta_{11},& \Delta_{12}  ,& \Delta_{22}
\end{array} \right]
\]

\begin{equation}
M_{2} = \left[ \begin{array}{rrrrrrrr}
x_{02} & 0 & x_{01} & - 0 & 0 & x_{00} & 0 & 0 \\
-x_{12} & x_{02} & -x_{11} & x_{01} & 0 & 0 & x_{00} & 0\\
x_{22} & 0 & x_{12} & x_{02} & x_{01} & x_{02} & 0 & x_{00}\\
0 & -x_{12} & 0 & -x_{11} & 0 & -x_{11} & -x_{01} & 0\\
0 & x_{22} & 0 & 0& -x_{11} & x_{12} & x_{02} & -x_{01}\\
0 & 0 & 0 & x_{22} & x_{12} & 0 & 0 & x_{02}
\end{array} \right]
\end{equation}

\begin{equation}
M_{3} = \left[ \begin{array} {rrr}
x_{01} & x_{00} & 0 \\
-x_{11} & -x_{01} & 0 \\
-x_{02} & 0 & x_{00} \\
x_{12} & x_{02} &  0 \\
-x_{22} & 0 & x_{02} \\
0 & -x_{02} & -x_{01}\\
0 & x_{12} & x_{11}\\
0 & -x_{22} & -x_{12}
\end{array} \right]
\end{equation}
\\ \\

\begin{center}
\textbf {3.  Resolutions of plane curves in the Veronese embedding.}
\end{center}

Let $C$ be a smooth(or irreducible) curve such that, $C \hookrightarrow
\mathbb{P}^{2}$. Hence C is given by a irreducible polynomial in three
variables. Now recall that $\mathbb{P}^{2} \hookrightarrow \mathbb{P}^{5}$.
Hence we have $C \hookrightarrow \mathbb{P}^{2} \stackrel{\sigma}
{\rightarrow } \mathbb{P}^{5}$. We will compute the syzygies of the
homogeneous ideal $\mathcal{I}_{\sigma(C)}$ using this embedding and the
resolution of the Veronese embedding above. Let $C$ be defined by the
polynomial $f$ of degree $d$ in three variables. Hence,
\[C = \mathcal{Z}(f(x_{0},x_{1},x_{2})) \]
Let,
\begin{displaymath}
f = \sum_{\substack{i+j+k = d}}a_{i,j,k} x_{0}^{i} x_{1}^{j} x_{2}^{k}
\end{displaymath}\\
\\

\textbf{3.1: Degree of $f$ is even}\\
\\

We have $d$ is even(say $2m$), and
\begin{displaymath}
f = \sum_{\substack{i+j+k = 2m}}a_{i,j,k} x_{0}^{i} x_{1}^{j} x_{2}^{k}
\end{displaymath}

\textbf{Lemma 2}: $Im(\theta)$ is a subalgebra of $K[x_{0},x_{1},x_{2}]$ and is generated by even polynomials.\\
Proof:\\ 
To prove thatl $f \in Im(\theta)$. We split $f$ in four parts, depending on
the parities of $i$,$j$,$k$, i.e., $f = f^{I} + f^{II} + f^{III} + f^{IV}$ with;
\[f^{I} = \sum_{\substack{i+j+k=d,\\ i,j,k \textrm{ even } }}a_{i,j,k} x_{0}^{i} x_{1}^{j} x_{2}^{k}\]
and so on.
 \\
\\
\emph{Case I} : When $i$, $j$, $k$ are all even, consider
\begin{displaymath}
F^{I} =  \sum_{\substack{i+j+k=d \\ i,j,k \textrm{ even }}} a_{ijk}x_{00}^{\frac{i}{2}}x_{11}^{\frac{j}{2}}x_{22}^{\frac{k}{2}}
\end{displaymath}

Notice that $\mathbf{\theta(F^{I}) = f^{I}}$
\\
\\
\emph{Case II}: When $i$ is even, $j$ and $k$ odd, consider
\begin{displaymath}
F^{II} = \sum_{\substack{i+j+k = d \\ i \textrm{ even } \\j,k \textrm{ odd } }} a_{ijk}x_{00}^{\frac{i}{2}}x_{11}^{\frac{j-1}{2}}x_{22}^{\frac{k-1}{2}}
x_{12}
\end{displaymath}
Similarly as Case I, $\mathbf{\theta(F^{II}) = f^{II}}$
\\
\\
\emph{Case III}: With $i$ is odd, $j$ is even, $k$ is odd consider
\begin{displaymath}
F^{III} = \sum_{\substack{i+j+k = d \\ j \textrm{ even } \\i,k \textrm{ odd } }} a_{ijk}x_{00}^{\frac{i-1}{2}}x_{11}^{\frac{j}{2}}x_{22}^{\frac{k-1}{2}}
x_{02}
\end{displaymath}
$ \mathbf{\theta(F^{III}) = f^{III}}$
\\
\\
\emph{Case IV}: $i$ is odd, $j$ is odd, $k$ is even consider,
\begin{displaymath}
F^{IV} = \sum_{\substack {i+j+k = d \\ k \textrm{ even } \\i,j \textrm{ odd } }}a_{ijk}x_{00}^{\frac{i-1}{2}}x_{11}^{\frac{j-1}{2}}x_{22}^{\frac{k}{2}}
x_{01}
\end{displaymath}
$\mathbf{\theta(F^{IV}) = f^{IV}}$
\\

Now let,
\begin{displaymath}
\mathbf{
F = F^{I} + F^{II} + F^{III} + F^{IV}}
\end{displaymath}
Then \begin{displaymath}
                \mathbf{\theta(F) = f }
                       \end{displaymath}
\\
\\

Hence $ \mathbf{f \in Im(\theta)}$.\\


\textbf{Lemma 3}: Let $G \in S \textrm { such that, } G \textrm{ homogeneous 
and  }$ $ \mathcal{Z}(\theta(F)) \subset \mathcal{Z}(\theta(G)) \subset \mathbb{P}^{2} $, where $F$ is a irreducible polynomial of even degree.
Then $G \in < F, \Delta_{i,j}: 0 \leq i \leq j \leq 2>$.\\
\\
\textbf{Proof}: Let $\theta(G) = g$, then $g$ is a homogeneous polynomial and,
\begin{displaymath}
\mathcal{Z}(f) \subset \mathcal{Z}(g)
\end{displaymath}
\begin{displaymath}
\Rightarrow g \in (f) \textrm{ as $C$ is a irreducible curve and hence $f$ is irreducible hence,}
\end{displaymath}
\begin{displaymath}
g= f.h \textrm{ for some $h$ homogeneous in } K[x_{0},x_{1},x_{2}]
\end{displaymath}
Now $f$ and $g$ are even degree implies that $h$ is of even degree hence, $ \exists H \in S$, homogeneous such that
$\theta(H) = h$. \\
Thus $\theta(G) = \theta(F). \theta(H)= \theta(F.H),$
\begin{displaymath}
\Rightarrow \theta(G-F.H) = 0
\end{displaymath}
\begin{displaymath}
\Rightarrow G-F.H \in ker(\theta)
\end{displaymath}
\begin{displaymath}
\Rightarrow G-F.H = \sum_{\substack{0\leq i \leq j \leq 2} }\Delta_{ij} S_{ij}
\textrm{  for some  }    S_{ij} \in S, S_{ij} \textrm{  homogeneous  }
\end{displaymath}
\begin{displaymath}
\Rightarrow G \in <   F, \Delta_{ij} : 0 \leq i \leq j \leq 2  >
\end{displaymath}
This completes the proof of the lemma.\\

Now recall $M_{2}$ and $M_{3}$ from equations (2) and (3), from now we will
denote them as below:
Let us denote the $i^{\textrm{th}}$ row of $M_{2}$ as $W_{i}$ and $j^{\textrm{th}}$ row of $M_{3}$ as $G_{j}$, for $1 \leq i \leq 8$ and $j=1,2,3$. Hence we get,
\[
M_{2} = \left[ \begin{array} {cccccccc}
W_{1},& W_{2},& W_{3},& W_{4},& W_{5},& W_{6},& W_{7},& W_{8} \end{array} \right] 
\hspace{1in} (*)
\]
\[
M_{3} = \left[ \begin{array} {ccc} G_{1},& G_{2},& G_{3}
\end{array} \right] \hspace{3in} (**)
\]
\\
\\
\\
\textbf{Theorem 4}: Let $C$ be an irreducible curve of even degree say $d=2m$, $m \geq 1$. The homogeneous ideal $\mathcal{I}_{C}$ of $\sigma(C)$  in $\mathbb{P}^{5}$ has the following resolution.\\
\begin{equation}
\begin{split}
0 \rightarrow & S(-m-4)^{\oplus 3}  \stackrel{M'_{4}}{\rightarrow}
S(-4)^{\oplus 3} \oplus S(-m-3)^{\oplus 8} \stackrel{M'_{3}}{\rightarrow}\\
& \stackrel{M'_{3}}{\rightarrow}
S(-3)^{\oplus 8} \oplus S(-m-2)^{\oplus 6} \stackrel{M'_{2}}{\rightarrow}
S(-2)^{\oplus 6} \oplus S(-m)\stackrel{M'_{1}} {\rightarrow}
S \rightarrow S/\mathcal{I}_{C} \rightarrow 0
\end{split}
\end{equation}
where,\\
\begin{equation}
M'_{1} = \left[ \begin{array} {rrrrrrrrrrr}
                              \left[  M_{1} \right]  ,& F
                       \end{array} \right]
\end{equation}
\\
\\
Also let,
\[N_{2} = \left[ \begin{array}{rrrrrrr}
           -F & 0 & 0  & 0  & 0 & 0 & \Delta_{00}\\
           0 & -F & 0  & 0  & 0 & 0 & \Delta_{01}\\
           0 & 0 & -F  & 0  & 0 & 0 & \Delta_{02}\\
           0 & 0 & 0  & -F  & 0 & 0 & \Delta_{11}\\
           0 & 0 & 0  & 0  & -F & 0 & \Delta_{12}\\
           0 & 0 & 0  & 0  & 0 & -F & \Delta_{22}\\

                     \end{array}\right]
\]
Let
\[
 N_{2} =\left[ \begin{array}{rrrrrrrrrrrrrrrrrrrrr}
                     U_{00},& U_{01},& U_{02},& U_{11},& U_{12},& U_{22}
                       \end{array} \right] ^ {T}
\]
\\
\begin{equation}
M'_{2} =\left[ \begin{array}{rrrrrrrrrrrrrrrrrrrrr}
                     W_{1}',& W_{2}',& W_{3}',& W_{4}',& W_{5}',& W_{6}',& W_{7}',& W_{8}',& U_{00},& U_{01},& U_{02},& U_{11},& U_{12},& U_{22}
                       \end{array} \right]
\end{equation}
\\
where
\[ W'_{i} = \left[ \begin{array}{cccccc}
                       W_{i} \\ 0
                       \end{array} \right]    \hspace{.5in} \forall i = 1, \ldots, 8
\]
with $W_{i}$ as in (*)
\\
Hence, 
\begin{center}
$M'_{2} = \left[  \begin{array} {ccccc} M_{2} & -FI_{6} \\ 0 & M_{1} \end{array} \right]$
\end{center}
Let
\[H_{i} = \left[ \begin{array}{cccccc} \left[F.I^{8}_{i}\right] \\ \left[ W_{i} \right] \end{array} \right] \]

where \[I^{k}_{i} = \left [ \begin{array} {rrrrrrrrr} 0,& 0,& \ldots ,& \stackrel { \textrm {$i^{th}$ position}} {1},& 0,& \ldots ,& 0 \end{array} \right]^{T} \textrm{is a } k \times 1 \textrm{ vector }
\]
\begin{equation}
M'_{3} = \left[ \begin{array} {rrrrrrrrrrrrr} G'_{1},& G'_{2},& G'_{3},& H_{1},& \ldots,&
H_{8} \end{array} \right]
\end{equation}
where,
\[
G'_{i} = \left[ \begin{array} {rrrrrrrrrrrrrrr} G_{i} \\ \left[ \bar{0} \right] \end{array} \right]  \hspace{1in} \textrm{ for } i=1,2,3\]\\
where $G_{i}$ as in (**) and $\left[ \bar{0} \right]$ is a $0$ matrix of appropriate dimension.\\
hence we have, 
\begin{center}
$M'_{3} = \left[  \begin{array} {ccccc} M_{3} & -FI_{8} \\ 0 & M_{2} \end{array} \right]$
\end{center}
Now let
\begin{equation}
M'_{4} = \left[  \begin{array} {ccccc} \left( \begin{array}  {cccccccccccccc}
\left[-F.I^{3}_{1} \right] \\ \left[G_{1}\right] \end{array}\right),&
\left( \begin{array}  {cccccccccccccc}
\left[-F.I^{3}_{2} \right] \\ \left[G_{2}\right] \end{array}\right),&
\left( \begin{array}  {cccccccccccccc}
\left[-F.I^{3}_{3} \right] \\ \left[G_{3}\right] \end{array}\right)
\end{array}\right]
\end{equation}
Hence we can write that,\\ 
\begin{center}
$M'_{4} = \left[  \begin{array} {ccccc} 
\left[-F.I^{3}\right] \\ \left[M_{3}\right] 
\end{array}\right]$
\end{center}
\textbf{Proof}:
\\
From Lemma 3, it is clear that \[ M_{1} = \left[ \begin{array}{rrrrrrrrrrr} \Delta_{00},& \Delta_{01},& \Delta_{02},& \Delta_{11},& \Delta_{12} ,& \Delta_{22},& F \end{array} \right] \]
\\
Now consider,
\[ \hspace{.130in} A= \left[ \begin{array} {rrrrrrrrrr} a_{00},& a_{01},& a_{02},& a_{11},& a_{12},& a_{22} \end{array} \right] \]
where  $a_{ij}  \in S$, homogeneous.
\[\textrm{And} \hspace{.125in}  B\in S, \textrm{ homogeneous }\]
 such that
\[\sum_{\substack{i,j}} a_{ij}.\Delta_{ij} +  B .F = 0\]
\[\Rightarrow \theta(B.F) = 0\]
\[\Rightarrow \theta(B).f = 0\]
\[\Rightarrow B \in < \Delta_{ij}: 0\leq i \leq j \leq 2> \]
Hence, $B = \sum (b_{ij} \Delta_{ij})$ for some homogeneous polynomials $b_{ij} \in S$.
\[ \Rightarrow \sum(a_{ij} + b_{ij}.F).\Delta_{ij} = 0 \]
Now if $a_{ij} + b_{ij}.F =0$ for all $(a_{ij},b_{ij})$ then such a $[ A, B]$ is generated by $U_{ij}$. If not then,\\
$\Rightarrow  \sum(a_{ij} + b_{ij}F) \in \textrm{Syz}^{1}( <  \Delta_{ij} : 0 \leq i \leq j \leq 2  > )$\\
Hence, the relations between $\Delta_{ij}$ and $F$ are generated by $U_{ij}:0 \leq i \leq j \leq 2$ and $W'_{k}:k=1,\ldots, 8.$\\

Hence we get,
\[ M'_{2} =\left[ \begin{array}{rrrrrrrrrrrrrrrrrrrrrrrrr}
                     W_{1}',& W_{2}',& W_{3}',& W_{4}',& W_{5}',& W_{6}',& W_{7}',& W_{8}',& U_{00},& U_{01},& U_{02},& U_{11},& U_{12},& U_{22}
                       \end{array} \right] \]
\\
Now consider \\
$A= \left[ \begin{array} {rrrrrrrrrrrr} a_{00},& a_{01},& a_{02},& a_{11},& a_{12},& a_{22}\end{array} \right]^{T} ,  a_{ij}  \in S$, $a_{ij}$ homogeneous $\forall 0 \leq i \leq j \leq 2$ and,
\[B = \left[ \begin{array} {rrrrrrrrrrrr} (b_{k}) \end{array} \right],  b_{k}\in S, \textrm{ homogeneous } \]
such that
\[ \sum_{\substack{0 \leq i \leq j \leq 2}} a_{ij}.U_{ij} + \sum_{\substack{1 \leq k \leq 8 }} b_{k} .W'_{k} = 0\]
\[ \Rightarrow \sum_{\substack{i,j} } a_{ij} \Delta_{ij}  = 0\]
As the last column of each $W'_{k}$, $k=1, \ldots, 8$ is zero and the last column of $U_{ij}$ is $\Delta_{ij}$ for $0\leq i\leq j \leq 2$
\[ \Rightarrow A \in < W_{k}: k=1, \ldots, 8 > \]
Let $A= \sum_{\substack{k}} (c_{k}W_{k})$, for some homogeneous polynomial, $c_{k} \in S$
\[ \Rightarrow -\sum_{\substack{k}}c_{k}W_{k}F.Id_{6}+ \sum_{\substack{k}}b_{k}W_{k} = 0\]
where $Id_{n}$ is a $n \times n$ identity matrix.
\[ \Rightarrow \sum_{\substack{i,k}} W_{k}(-c_{k}F + b_{k}) =0\]

Hence if $-c_{k}.F + b_{k} = 0$ for all $k$, this implies $b_{k} = c_{k}.F$ for all $k$ then such $(b_{k},a_{ij})$ are generated by $< \left[ \left[ F.[I^{8}_{i}]\right] , \left[ W_{i} \right] \right]>$ for $i=1,\ldots,8$.
And if not then,
$\Rightarrow {\left[ (-c_{k}F+b_{k})I_{k} \right]}_{k=1, \ldots ,6} \in {\textrm{Syz}}^{1}(<W_{j}:j=1,\ldots,8>)$.\\

Hence the relations between $W'_{k}$ and $U_{ij}$ are generated by $G'_{i}$ and $H_{k}$. Hence we get,
\[
M'_{3} = \left[ \begin{array} {rrrrrrrrrrrrr} G'_{1},& G'_{2},& G'_{3},&  H_{1},& \ldots,& H_{8}
                      \end{array} \right]
\]

Now consider \\

$ A= \left[ \begin{array} {rrrrrrrrrrrr} a_{1},& a_{2},& a_{3},& a_{4},& a_{5},& a_{6},& a_{7},& a_{8} \end{array} \right]^{T},  a_{i}  \in S$, homogeneous for $i=1, \ldots, 8$\\

$B = \left[ \begin{array} {rrrrrrrrrrrr} (b_{k}) \end{array} \right]$  $b_{k}\in S$, homogeneous for $k=1,2,3$ such that
\[\sum_{\substack{i}} a_{i}.H_{i} + \sum_{\substack{k}}b_{k} .G'_{k} = 0\]
\[ \Rightarrow \sum_{\substack{i} } a_{i} W_{i}  = 0\]
As the last six columns of each $G'_{k}$, $k=1,2,3$ are zero.
\[\Rightarrow A \in < G_{p}: p=1,2,3 > \]
Let $A= \sum_{\substack{k}} (c_{p}G_{p})$, for some homogeneous polynomial, $c_{k} \in S$.\\
Then we have, $\sum_{\substack{p}}(c_{p}G_{p}).(F.Id_{8}) + \sum_{\substack{k}}b_{k}.G_{k} =0$
\[\Rightarrow \sum_{\substack{p}}\left(c_{p}.F.Id_{8} + b_{p} \right) G_{p} = 0\]
Now if  $c_{p}.F + b_{p} = 0$ for every $p$, then $ b_{p} = -c_{p}.F$ for all $p$, then we can say that
$([b_{p}],[c_{p}])$  is generated by  $< \left( \left[-F.I^{3}_{i}\right], \left[I^{3}_{i}\right] \right)   :   i=1,2,3>$, hence
$([b_{p}],[a_{k}])$ is generated by $< \left( \left[-F.I^{3}_{i}\right],\left[G_{i}\right]\right) i=1,2,3>$\\

Also from theorem 1 we have that ${G'_{k}:k=0,1,2}$ are independent. Hence \\
$\textrm{Syz}^{1}(<G'_{i}, H_{j} : i=1,2,3 \textrm{ and } j=1,\ldots 8>) =< \left( \left[-F.I^{3}_{i}\right],\left[G_{i}\right]\right)  :i=1,2,3>$\\
Hence,
\[M'_{4} = \left[ \left( \begin{array} {cccccccccccccc}
                   \left[-F.I^{3}_{i}\right]
                  \\ \left[G_{i}\right]

                               \end{array} \right) \right]\]
$i=1,2,3$
\\

\textbf{3.2: Degree of f is odd}\\
\\

Recall
\begin{displaymath}
f = \sum_{\substack{i+j+k = d}}a_{i,j,k} x_{0}^{i} x_{1}^{j} x_{2}^{k}
\end{displaymath}
\\
Now let $f_{0}= x_{0}.f$    ,   $f_{1}= x_{1}.f$   ,    $f_{2}= x_{2}.f$\\
Then $f_{n}$ is of even degree and hence according to Case A,
  $f_{n} \in Im(\theta)$  for  $n=0,1,2$  \\

\textbf{Lemma 5}: $\hspace{.1in}$    Let $G \in k[x_{00},x_{01},x_{02},x_{11},x_{12},x_{22}] $ such that, $G$
homogeneous and \\
$\mathcal{Z}(\theta(F_{0})) \cap \mathcal{Z}(\theta(F_{1})) \cap \mathcal{Z}(\theta(F_{2})) \subset \mathcal{Z}( \theta(G) ) \subset
\mathbb{P}^{2} $.
Then $G \in <  F_{k}, \Delta_{i,j}:0 \leq k \leq 2, 0 \leq i \leq j \leq 2   >$.\\
\\
\textbf{Proof}:$\hspace{.1in}$ Now let $\theta(G) = g$, then degree($g$) is even.
\\
\\
\begin{displaymath}
\mathcal{Z}(f_{0})\cap \mathcal{Z}(f_{1})\cap \mathcal{Z}(f_{2})\subset \mathcal{Z}(g)
\end{displaymath}
\\
\begin{displaymath}
\Rightarrow \mathcal{Z}(f)\subset \mathcal{Z}(g)
\end{displaymath}
\\
\begin{displaymath}
\Rightarrow g \in (f) \textrm{as $C$ is an irreducible curve and hence $f$ is irreducible}
\end{displaymath}
\\
\begin{displaymath}
\Rightarrow g= f.h \textrm{ for some } h \in  k[x_{0},x_{1},x_{2}]
\end{displaymath}
\\
\\
\[ \Rightarrow h \neq 1 \textrm{  as degree  } f \textrm{ is odd while degree } g  \textrm{ is even} \]
\\
\[\Rightarrow g = \sum_{\substack{i=0,1,2}}f_{i}h_{i} \textrm{ for some  homogeneous even degree polynomials } h_{i} \in k[x_{0},x_{1},x_{2}]\]
\\
\\
$\Rightarrow G = \sum_{\substack{i=0,1,2}}F_{i}H_{i},$  where  $\theta(H_{i}) = h_{i}$  $\forall i=0,1,2.$ \\
$\hspace{.5in}$ Such a $H_{i}$, exists as the degree of $h_{i}$ is even.
\\
\begin{displaymath}
\Rightarrow \theta(G - \sum_{\substack{i=0,1,2}}F_{i}H_{i}) = 0
\end{displaymath}
\\
\begin{displaymath}
\Rightarrow G-\sum_{\substack{i=0,1,2}}F_{i}H_{i} \in ker(\theta)
\end{displaymath}
\\
\begin{displaymath}
\Rightarrow G =\sum_{\substack{i=0,1,2} } F_{i}H_{i} + \sum_{\substack{i,j=0,1,2}} \Delta_{ij} S_{ij} \hspace{.1in} \textrm{for some} \hspace{.1in}S_{ij}  \in k[x_{00},\ldots ,x_{22}]
\end{displaymath}
\\
\begin{displaymath}
\Rightarrow G \in < F_{k}, \Delta_{ij} : i,j,k=0,1,2>
\end{displaymath}
\\
\\

\textbf{Lemma 6}: $Im(\theta)$ is a subalgebra of $K[x_{0},x_{1},x_{2}]$ and is generated by even polynomials.\\
Proof:\\
Like in the case of degree of $f$ being even we split $f$ in four parts 
depending on the parities of $i$,$j$,$k$.\\
Case I: $i$, $j$, $k$ are all odd. Let
\\
\begin{displaymath}
\textrm{ Let } {\mathit{h}}_{I} = \sum_{\substack{i,j,k}}a_{ijk} x_{00}^{\frac {i-1}{2}}x_{11}^{\frac {j-1}{2}}x_{22}^{\frac {k-1}{2}}
\end{displaymath}
\\
\begin{displaymath}
{F_{0}}^{I} = \sum_{\substack{i+j+k =d}}a_{i,j,k} x_{00}^{\frac{i+1}{2}}x_{11}^{\frac{j-1}{2}}x_{22}^{\frac{k-1}{2}}x_{12}
\end{displaymath}
\\
\begin{displaymath}
{F_{1}}^{I} = \sum_{\substack{i+j+k =d}}a_{i,j,k} x_{00}^{\frac{i-1}{2}}x_{11}^{\frac{j+1}{2}}x_{22}^{\frac{k-1}{2}}x_{02}
\end{displaymath}
\\
\begin{displaymath}
{F_{2}}^{I} = \sum_{\substack{i+j+k =d}}a_{i,j,k} x_{00}^{\frac{i-1}{2}}x_{11}^{\frac{j-1}{2}}x_{22}^{\frac{k+1}{2}}x_{01}
\end{displaymath}
\\
Then,
\begin{displaymath}
{F_{0}}^{I} = x_{00}x_{12}{\mathit{h}}_{I}
\end{displaymath}
\\
\begin{displaymath}
{F_{1}}^{I} = x_{11}x_{02}{\mathit{h}}_{I}
\end{displaymath}
\\
\begin{displaymath}
{F_{2}}^{I} = x_{22}x_{01}{\mathit{h}}_{I}
\end{displaymath}
\\
Case II: $i$ odd, $j$ even, $k$ even. Now
\begin{displaymath}
\textrm{ Let } {\mathit{h}}_{II} = \sum_{\substack{i,j,k }}a_{ijk} x_{00}^{\frac {i-1}{2}}x_{11}^{\frac {j}{2}}x_{22}^{\frac {k}{2}}
\end{displaymath}
\\
\begin{displaymath}
{F_{0}}^{II} = \sum_{\substack{i+j+k =d}}a_{i,j,k} x_{00}^{\frac{i+1}{2}}x_{11}^{\frac{j}{2}}x_{22}^{\frac{k}{2}}
\end{displaymath}
\\
\begin{displaymath}
{F_{1}}^{II} = \sum_{\substack{i+j+k =d}}a_{i,j,k} x_{00}^{\frac{i-1}{2}}x_{11}^{\frac{j}{2}}x_{22}^{\frac{k}{2}}x_{01}
\end{displaymath}
\\
\begin{displaymath}
{F_{2}}^{II} = \sum_{\substack{i+j+k =d}}a_{i,j,k} x_{00}^{\frac{i-1}{2}}x_{11}^{\frac{j}{2}}x_{22}^{\frac{k}{2}}x_{02}
\end{displaymath}
\\
Then,
\begin{displaymath}
{F_{0}}^{II} = x_{00}{\mathit{h}}_{II}
\end{displaymath}
\\
\begin{displaymath}
{F_{1}}^{II} = x_{01}{\mathit{h}}_{II}
\end{displaymath}
\\
\begin{displaymath}
{F_{2}}^{II} = x_{02}{\mathit{h}}_{II}
\end{displaymath}
\\
Case III: $i$ even, $j$ odd, $k$ even. Now
\\
\begin{displaymath}
\textrm{ Let } {\mathit{h}}_{III} = \sum_{\substack{i,j,k}}a_{ijk} x_{00}^{\frac {i}{2}}x_{11}^{\frac {j-1}{2}}x_{22}^{\frac {k}{2}}
\end{displaymath}
\\
\begin{displaymath}
{F_{0}}^{III} = \sum_{\substack{i+j+k =d}}a_{i,j,k} x_{00}^{\frac{i}{2}}x_{11}^{\frac{j-1}{2}}x_{22}^{\frac{k}{2}}x_{01}
\end{displaymath}
\\
\begin{displaymath}
{F_{1}}^{III} = \sum_{\substack{i+j+k =d}}a_{i,j,k} x_{00}^{\frac{i}{2}}x_{11}^{\frac{j+1}{2}}x_{22}^{\frac{k}{2}}
\end{displaymath}
\\
\begin{displaymath}
{F_{2}}^{III} = \sum_{\substack{i+j+k =d}}a_{i,j,k} x_{00}^{\frac{i}{2}}x_{11}^{\frac{j-1}{2}}x_{22}^{\frac{k}{2}}x_{12}
\end{displaymath}
\\
Then,\\
\begin{displaymath}
{F_{0}}^{III} = x_{01}{\mathit{h}}_{III}
\end{displaymath}
\\
\begin{displaymath}
{F_{1}}^{III} = x_{11}{\mathit{h}}_{III}
\end{displaymath}
\\
\begin{displaymath}
{F_{2}}^{III} = x_{12}{\mathit{h}}_{III}
\end{displaymath}
\\
Case IV: $i$ even, $j$ even, $k$ odd. Now\\
\\
\begin{displaymath}
\textrm{ Let } {\mathit{h}}_{IV} = \sum_{\substack{i,j,k}}a_{ijk} x_{00}^{\frac {i}{2}}x_{11}^{\frac {j}{2}}x_{22}^{\frac {k-1}{2}}
\end{displaymath}
\\
\begin{displaymath}
{F_{0}}^{IV} = \sum_{\substack{i+j+k =d}}a_{i,j,k} x_{00}^{\frac{i}{2}}x_{11}^{\frac{j}{2}}x_{22}^{\frac{k-1}{2}}x_{02}
\end{displaymath}
\\
\begin{displaymath}
{F_{1}}^{IV} = \sum_{\substack{i+j+k =d}}a_{i,j,k} x_{00}^{\frac{i}{2}}x_{11}^{\frac{j}{2}}x_{22}^{\frac{k-1}{2}}x_{12}
\end{displaymath}
\\
\begin{displaymath}
{F_{2}}^{IV} = \sum_{\substack{i+j+k =d}}a_{i,j,k} x_{00}^{\frac{i}{2}}x_{11}^{\frac{j}{2}}x_{22}^{\frac{k+1}{2}}
\end{displaymath}
\\
Then,
\begin{displaymath}
{F_{0}}^{IV} = x_{02}{\mathit{h}}_{IV}
\end{displaymath}
\\
\begin{displaymath}
{F_{1}}^{IV} = x_{12}{\mathit{h}}_{IV}
\end{displaymath}
\\
\begin{displaymath}
{F_{2}}^{IV} = x_{22}{\mathit{h}}_{IV}
\end{displaymath}
\\
\\
\\
$F_{n} = {F_{n}}^{I} + {F_{n}}^{II} + {F_{n}}^{III} + {F_{n}}^{IV}$
$ \forall n =0,1,2 $ \\
Also notice $\theta(F_{n}) = f_{n}$  for  $n=0,1,2$ \\
\\
\textbf{Theorem 7}: Let $C$ be an irreducible curve of odd degree say $d = 2m - 1$, for $m \geq 2$. The ideal $\mathcal{I}_{C}$ of $\sigma(C)$  in $\mathbb{P}^{5}$ has the following resolution.
\begin{equation}
\begin{split}
0 \rightarrow & S(-m-4)  \stackrel{\beta_{4}}{\rightarrow}
S(-4)^{\oplus 3} \oplus  S(-m-2)^{\oplus 6} \stackrel{\beta_{3}}{\rightarrow}\\
& \stackrel{\beta_{3}}{\rightarrow}
S(-3)^{\oplus 8} \oplus S(-m-1)^{\oplus 8} \stackrel{\beta_{2}}{\rightarrow}
S(-2)^{\oplus 6} \oplus  S^{\oplus 3}(-m) \stackrel{\beta_{1}} {\rightarrow} S \rightarrow S/\mathcal{I}_{C} \rightarrow 0
\end{split}
\end{equation}
\textbf{Proof}:
\\
From Lemma 3 and Lemma 5, it is clear that 
\[ \beta_{1} = \left[ \begin{array}{rrrrrrrrrrrrrrrrrr} \Delta_{00},& \Delta_{01},& \Delta_{02},& \Delta_{11},& \Delta_{12},& \Delta_{22},& F_{0},& F_{1},& F_{2}, \end{array} \right] \]
\\
Now consider 
$A= \left[ \begin{array} {rrrrrrrrrrrrrrrr} a_{00},& a_{01},& a_{02},& a_{11},& a_{12},& a_{22} \end{array} \right]$,  $a_{ij}  \in S$, homogeneous  $\forall 0 \leq i \leq \j \leq 2$ and $b= \left[ \begin{array} {rrrrrrrrrrrrrrrr} b_{0},& b_{1},& b_{2} \end{array} \right]$ where $b_{l} \in S$, homogeneous, for $k=0,1,2$ such that,   
\begin{equation}
\sum_{\substack{i,j}} a_{ij}.\Delta_{ij} +  \sum_{k}b_{k} .F_{k} = 0
\end{equation}
\[\Rightarrow \theta(\sum_{k} (b_{k}.F_{k})) = 0\]
\[\Rightarrow \sum_{k}(\theta(b_{k}).f_{k})= 0\]
\[\Rightarrow \sum_{k}(\theta(b_{k}).f.x_{k})= 0\]
\[\Rightarrow \sum_{k}(\theta(b_{k}).x_{k})= 0\]
Let $\theta(b_{k}) = B_{k}$, then degree of $B_{k}$ is even.
Then,
\[ B=(B_{0},B_{1},B_{2})^{T} \in \textrm{Syz}^{1} (x_{0},x_{1},x_{2})\]
Now by simple computation we get 
\begin{displaymath}
\textrm{Syz}^{1} (x_{0},x_{1},x_{2}) = <  \left( \begin{array} {rrrrrrrrrrrr} x_{1} \\ -x_{0} \\ 0  \end{array} \right) , 
\left( \begin{array} {rrrrrrrrrrrr} x_{2} \\ 0 \\ -x_{0} \end{array} \right), 
\left( \begin{array} {rrrrrrrrrrrr} 0 \\ x_{2} \\ -x_{1} \end{array} \right) > 
\end{displaymath}
hence $B \in < Y_{0},Y_{1},Y_{2}> $
where,
\[ \left.\begin{array}{lllll}
 Y_{0} =  \left( \begin{array} {rrrrrrrrrrrr} x_{1},& -x_{0},& 0  \end{array} \right) \\
\\ Y_{1} =  \left( \begin{array} {rrrrrrrrrrrr} x_{2},& 0,& -x_{0} \end{array} \right) \\
\\ Y_{2} =  \left( \begin{array} {rrrrrrrrrrrr} 0,& x_{2},& -x_{1} \end{array} \right) 
\end{array} \right. \]
But degree of $B_{k}$ is even, hence,
$B \in < x_{k}Y_{l} : k,l =0,1,2 >$.\\
Hence, $ ( b_{0},b_{1},b_{2} ) \in < Y_{lk} : k,l =0,1,2 >$\\
\\
where 
\[Y_{00} =  \left( \begin{array} {rrrrrrrrrrrr} x_{01} ,& -x_{00} ,& 0  \end{array} \right)  \]
\[Y_{01} =  \left( \begin{array} {rrrrrrrrrrrr} x_{11} ,& -x_{01} ,& 0  \end{array} \right)  \]
\[Y_{02} =  \left( \begin{array} {rrrrrrrrrrrr} x_{12} ,& -x_{02} ,& 0  \end{array} \right)  \]

\[Y_{10} =  \left( \begin{array} {rrrrrrrrrrrr} x_{02} ,& 0 ,& -x_{00}   \end{array} \right)  \]
\[Y_{11} =  \left( \begin{array} {rrrrrrrrrrrr} x_{12} ,& 0 ,& -x_{01}   \end{array} \right)  \]
\[Y_{12} =  \left( \begin{array} {rrrrrrrrrrrr} x_{22} ,& 0 ,& -x_{02}   \end{array} \right)  \]

\[Y_{20} =  \left( \begin{array} {rrrrrrrrrrrr} 0 ,&  x_{02} ,& -x_{01}  \end{array} \right)  \]
\[Y_{21} =  \left( \begin{array} {rrrrrrrrrrrr} 0 ,&  x_{12} ,& -x_{11}  \end{array} \right)  \]
\[Y_{22} =  \left( \begin{array} {rrrrrrrrrrrr} 0 ,&  x_{22} ,& -x_{12}  \end{array} \right)  \]
Also note that,
\[Y_{02}=Y_{11} - Y_{20}\]
Now substituting all $Y_{ij}$ for $i,j=0,1,2$ except for $Y_{02}$ for $b$ in (10) we get, the following 8 vectors,\\
\\
\[
\left. \begin{array}{lllllllllllllllllllllllllllllllllllllllll}
V_{1} = \left[ \begin{array}{cccccccccccccccccccc} 0,& 0,& -x_{00}{\mathit{h}}_{I},& 
0,& {\mathit{h}}_{IV},& {\mathit{h}}_{III},& \left[ Y_{00}\right] \end{array} 
\right]^T \\ \\
V_{2} = \left[  \begin{array}{cccccccccccccccccccc}0,& 0,&{\mathit{h}}_{IV},& 0 
,&-x_{11}{\mathit{h}}_{I} ,& -{\mathit{h}}_{II} ,& \left[ Y_{01}\right] \end{array}\right]^T \\ \\
V_{3} = \left[ \begin{array} {ccccccccccccccccc}0 ,& x_{00}{\mathit{h}}_{I} ,& 0 ,& {\mathit{h}}_{IV} ,& {\mathit{h}}_{III} ,& 0 ,&\left[ Y_{10}\right] \end{array}\right]^T \\ \\   
V_{4} = \left[\begin{array}{ccccccccccccccccc} x_{00}{\mathit{h}}_{I} ,& {\mathit{h}}_{IV} ,& 0 ,& 0 ,&-{\mathit{h}}_{II} ,& -x_{22}{\mathit{h}}_{I} ,& \left[ Y_{11}\right] \end{array}\right]^T \\ \\
V_{5} = \left[ \begin{array}{cccccccccccccccccccccc}
 0 ,& -{\mathit{h}}_{III} ,& 0 ,& -{\mathit{h}}_{II} ,& -x_{22}{\mathit{h}}_{I} ,& 0,& \left[ Y_{12}\right] \end{array} 
\right]^T \\ \\
V_{6} = \left[ \begin{array} {cccccccccccccccccccccc}0 ,& {\mathit{h}}_{IV} ,& {\mathit{h}}_{III} ,& x_{11}{\mathit{h}}_{I} ,& 0 ,& -x_{22}{\mathit{h}}_{I} ,& \left[ Y_{20}\right] \end{array}\right]^T \\ \\
V_{7} = \left[ \begin{array}{ccccccccccccccccc}{\mathit{h}}_{IV} ,& x_{11}{\mathit{h}}_{I} ,& -{\mathit{h}}_{II} ,& 0 ,& 0 ,& 0 ,& \left[ Y_{21}\right] \end{array}\right] ^T \\ \\
V_{8} = \left[ \begin{array}{cccccccccccccccccccccc}-{\mathit{h}}_{III},& -{\mathit{h}}_{II},& x_{22}\mathit{h}_{I},& 0,& 0,& 0,& \left[ Y_{22}\right] \end{array}\right] ^T \end{array} \right. \]

Let 
\[\beta'_{2} = \left[ \begin{array} {ccccccccc} \left[V_{1}\right],& \left[V_{2}\right],&\left[V_{3}\right],&\left[V_{4}\right],&\left[V_{5}\right],&\left[V_{6}\right],&\left[V_{7}\right],&\left[V_{8}\right] \end{array}\right] \]
\\
Now all the relations between $F_{n}$'s and $\Delta_{ij}$'s are generated by $V_{k}$'s and $W'_{l}$'s and all the relations between only $\Delta_{ij}$'s are generated by $W_{l}$'s. Hence all relations between ${ {F_{n}}, {\Delta_{jk} } }$ are generated by ${V_{k},W'_{l}}$.\\ 
Hence $\textrm{Syz}^{1}(<F_{n},\Delta_{ij}>) = < V_{k},W'_{l}>$. 
Hence \[\beta_{2} = \left( \left[W'_{1}\right]  \left[W'_{2}\right] \ldots \left[W'_{8}\right] \left[V_{1} \right]\ldots \left[V_{8}\right] \right)\]
where $W'_{k} = \left[ \left[W_{k}\right] \left[ \bar{0} \right] \right]$ with $\left[ \bar{0} \right]$ a $1 \times 3$ zero vector\\
Hence we can write $\beta'_{2}$ as\\,
\[\beta_{2} = \left[ \begin{array} {ccccccccc} \left[M_{2}\right]& H \\ \bar{0} & Y \end{array}\right] \]

Now consider, $\bar{A} = \left(a_{i}\right)$ such that $a_{i} \in S$ homogeneous and $\bar{B} = \left(b_{k}\right)$, where $b_{k} \in S$, homogeneous such that
\begin{equation}
\sum_{\substack{i} } a_{i}V_{i} + \sum_{\substack{k}} b_{k}W'_{k} = 0
\end{equation}
Now as all the entries in the last 3 columns in each of $W_{i}'$ are zero we have, 
\[ \sum_{\substack{i} } A_{i}Y_{ij} = 0\].
\\ 
Now it can be computed that ${\textrm{Syz}}^{1}( Y_{ij}) = < L_{i} : 1 \leq i \leq 6>$ where,
\[
\left. \begin{array} {lllllllllllllllllllllllllllllllllll}
L_{1} =  \left[\begin{array}{cccccccc} x_{02} ,&  0  ,&  -x_{01} ,& 0 ,& 0 ,& x_{00} ,& 0 ,& 0 \end{array}\right]\\
\\ L_{2} =  \left[\begin{array}{cccccccc} x_{12} ,&  x_{02}  ,& -x_{11} ,& -x_{01} ,& 0 ,& x_{01} ,& x_{00} ,& 0 \end{array}\right]\\
\\L_{3} =  \left[\begin{array}{cccccccc} x_{22} ,&  0  ,& -x_{12} ,& x_{02} ,& -x_{01} ,& 0 ,& 0 ,& x_{00} \end{array}\right]\\
\\ L_{4} =  \left[\begin{array}{cccccccc} 0 ,&  x_{12} ,& 0 ,& -x_{11} ,& 0 ,& 0 ,& x_{01} ,& 0 \end{array}\right]\\
\\ L_{5} =  \left[\begin{array}{cccccccc} 0  ,& x_{22} ,& 0 ,& 0 ,& -x_{11} ,& -x_{12} ,& x_{02} ,& x_{01}  \end{array}\right]\\
\\ L_{6} =  \left[\begin{array}{cccccccc} 0 ,& 0 ,& 0 ,& x_{22} ,&  -x_{12} ,& -x_{22} ,& 0 ,& x_{02} \end{array}\right]
\end{array} \right.
\]
\\
So substituting $K_{i}$, $i=0, \ldots, 6$ for $\bar{B}$ in (11) we get the following 6 vectors,

\[
\left.
\begin{array}{llllllllllllllllllllllllllllll}
K_{1}=  \left[  \begin{array}{ccccccccc}   0 ,& 0 ,& 0 ,& x_{00}\mathit{h}_{I} ,& 0 ,& 0 ,& -\mathit{h}_{IV} ,& \mathit{h}_{III} ,& \left[ L_{1} \right]\end{array} \right] ^T \\
\\K_{2}=  \left[ \begin{array}{ccccccccc} 0 ,& 0 ,& x_{00}\mathit{h}_{I} ,& 0 ,& -\mathit{h}_{III} ,& -\mathit{h}_{IV} ,&  x_{11}\mathit{h}_{I} ,& \mathit{h}_{II} ,& \left[ L_{2} \right] \end{array} \right]^T \\
\\K_{3}=  \left[  \begin{array}{ccccccccc} -x_{00}\mathit{h}_{I} ,& -\mathit{h}_{IV} ,& 0
,&  -\mathit{h}_{III} ,& 0 ,& \mathit{h}_{III} ,&  \mathit{h}_{II} ,& x_{22}\mathit{h}_{I} ,& \left[ L_{3} \right] \end{array} \right] ^T\\
\\K_{4}=  \left[ \begin{array} {ccccccccc}  0 ,& 0 ,& -\mathit{h}_{IV} ,& -x_{11}\mathit{h}_{I} ,& \mathit{h}_{II} ,& x_{11}\mathit{h}_{I} ,& 0 ,& 0 ,& \left[ L_{4} \right] \end{array}\right] ^T\\
\\K_{5}=  \left[ \begin{array} {ccccccccc}  -\mathit{h}_{IV} ,&  -x_{11}\mathit{h}_{I} ,& \mathit{h}_{III} ,& \mathit{h}_{II} ,-x_{22}\mathit{h}_{I} ,& 0 ,& 0 ,& 0 & \left[ L_{5} \right] \end{array} \right] ^T\\
\\K_{6}=  \left[  \begin{array}{ccccccccc}  \mathit{h}_{III} ,& \mathit{h}_{II} ,&  0 ,&  0 ,&  0 ,& -x_{22}\mathit{h}_{I} ,&  0  ,& 0 ,& \left[ L_{6} \right] \end{array}\right] ^T
\end{array}
\right.\]
\\
\\
Now all the relations between $V_{i}$'s and $W'_{j}$'s are generated by $ \{ K_{l}$'s, $G'_{k}$'s, $1 \leq l \leq 6 , k=1,2,3\}$ and all the relations between only $W'_{j}$'s (which are actually $W_{j}$) are generated by $G'_{k}$'s. Hence we have that all relations between $\{ \{V_{i}\}, \{W'_{j} \} \}$ are generated by $\{ K_{l}$'s,$G'_{k}$'s $\}$. So ,\\
$\textrm{Syz}^{1}(<V_{i},W'_{j}>) = < K_{l},G'_{k}>$. So we get that,
\[\beta_{3} = \left[ \begin{array}{cccccccc}\left[ G'_{0}\right] & \left[G'_{1}\right] & \left[G'_{2}\right] &\left[K_{1}\right] & \ldots &\left[ K_{6}\right] \end{array}\right]\]
where, $G'_{i} = \left[ \begin{array} {cccccccc}\left[G_{i}\right]& \left [\bar{0}\right] \end{array}\right]$ where $\left[ \bar{0} \right]$ is an appropriate dimensional zero matrix. \\
Hence we can write that, 
\[\beta_{3} = \left[ \begin{array}{cccccccc}M_{3} & L \\ 0 & K \end{array}\right]\]
Now consider $\bar{A} = \left(A_{i}\right)$, such that $A_{i} \in S$, homogeneous and $\bar{B} = \left(B_{k}\right)$, such that  $B_{k} \in S$, homogeneous such that,
\begin{equation}
\sum_{\substack{l} } A_{l}K_{l} + \sum_{\substack{k}} B_{k}G'_{k} = 0
\end{equation}
Hence we have, \[ \sum_{\substack{l} } A_{l}K'^{T}_{l} = 0 \]
(as the last eight columns of $G'_{i}$'s are zero entries)
\\
Now it can be computed that ${\textrm{Syz}}^{1}( K'_{l}) = < J'>$ where,
\[ J' = \left[ \begin{array} {ccccccccccccccccc}
{x_{12}}^2 - x_{11}x_{22} \\
-x_{02}x_{12} + x_{01}x_{22} \\
x_{11}x_{02} - x_{01}x_{12} \\
 {x_{02}}^2-x_{00}x_{22} \\
-x_{01}x_{02}+x_{00}x_{12} \\
{x_{01}}^{2}-x_{00}x_{11} \end{array} \right] \]
\\
Like earlier, substituting $J'$ in (12) we get $J$.
\[J =  \left[ \begin{array} {ccccccccccc}
 -x_{00}x_{12}\mathit{h}_{I} -x_{00}\mathit{h}_{II} - x_{01}\mathit{h}_{III} - x_{02}\mathit{h}_{IV} \\
\\
-x_{11}x_{02}\mathit{h}_{I} + x_{01}\mathit{h}_{II} + x_{11}\mathit{h}_{III} + x_{12}\mathit{h}_{IV} \\
\\
-x_{01}x_{22}\mathit{h}_{I} -x_{02}\mathit{h}_{II} - x_{12}\mathit{h}_{III} - x_{22}\mathit{h}_{IV} \\
\\
 \left[ J'  \right]  \end{array} \right]   \]
Now all the relations between $K_{l}$'s and $G'_{k}$'s are generated by $J$ and there are no relations between only $G'_{k}$'s as there are no non-trivial relations between $G_{k}$'s. Hence all relations between ${ {K_{l}}, {G'_{k} } }$ are generated by $J$. Hence $\textrm{Syz}^{1}(<K_{l},G'_{k}>) = < J >$. Hence\\
\[\beta_{4} = \left[ J \right]\]\\
This completes the proof of the theorem.
\\
\\
\emph{Acknowledgments}: \\
I used the program Singular on some curves of degree less than 5 to get
an idea about the resolutions of curves in the general case.\\
I would like to thank Prof. D.S.Nagaraj, for suggesting this problem
and for his suggestions and comments and a careful reading of the
manuscript.\\
Thanks are due to Prof Clare D'Cruz  for explaining the use of the
program Singular.\\
I also thank  Prof  Tony J. Puthenpurakal for valuable comments and
suggestions.\\
I want to thank The Institute of Mathematical Sciences for supporting me when 
this problem was being addressed.

 \end{document}